\documentclass[twocolumn, amsthm]{autart}

\usepackage{graphicx}
\usepackage{epsfig}

\usepackage{amsmath}
\usepackage{amsfonts}
\usepackage{amsthm}
\usepackage{mathrsfs}
\usepackage{xspace}
\usepackage{xparse}
\usepackage[inline]{enumitem}
\usepackage{etoolbox}
\usepackage{mathtools}

\usepackage{algorithm}
\usepackage{algpseudocode}

\usepackage[space]{cite}

\usepackage{tabularx}

\usepackage{lipsum}
\usepackage{color, soul}

\usepackage{theoremref}

\begin{document}

\begin{frontmatter}

\title{State-Based Confidence Bounds for Data-Driven \\ Stochastic Reachability Using Hilbert Space Embeddings\thanksref{footnoteinfo}}

\thanks[footnoteinfo]{
Corresponding author Meeko Oishi (\texttt{oishi@unm.edu}).
This material is based upon work supported by the National Science Foundation under NSF Grant Number CNS-1836900.
Any opinions, findings, and conclusions or recommendations expressed in this
material are those of the authors and do not necessarily reflect the views
of the National Science Foundation.
This research was  
supported in part by the Laboratory Directed Research and Development program at Sandia National Laboratories, a multimission laboratory managed and operated by National Technology and Engineering Solutions of Sandia, LLC., a wholly owned subsidiary of Honeywell International, Inc., for the U.S. Department of Energy's National Nuclear Security Administration under contract DE-NA-0003525.
The views expressed in this article do not necessarily represent the views of the U.S. Department of Energy or the United States Government.
This work was also supported in part by the University of New Mexico's Vice President for Research Office.
}

\author[UNM]{Adam J. Thorpe}\ead{ajthor@unm.edu},
\author[UNM]{Kendric R. Ortiz}\ead{kendric@unm.edu},
\author[UNM]{Meeko M. K. Oishi}\ead{oishi@unm.edu}

\address[UNM]{Electrical and Computer Engineering, University of New Mexico, Albuquerque, NM, USA}


\begin{abstract}
In this paper, we compute finite sample bounds for data-driven approximations of the solution to stochastic reachability problems.
Our approach uses a nonparametric technique known as kernel distribution embeddings, and provides probabilistic assurances of safety for stochastic systems in a model-free manner. 
By implicitly embedding the stochastic kernel of a Markov control process in a reproducing kernel Hilbert space, we can approximate the safety probabilities for stochastic systems with arbitrary stochastic disturbances as simple matrix operations and inner products. 
We present finite sample bounds for point-based approximations of the safety probabilities through construction of probabilistic confidence bounds that are state- and input-dependent. One advantage of this approach is that the bounds are responsive to non-uniformly sampled data, meaning that tighter bounds are feasible in regions of the state- and input-space with more observations.
We numerically evaluate the approach, and demonstrate its efficacy on a neural network-controlled pendulum system.
\end{abstract}

\end{frontmatter}


\section{Introduction}

In expensive, high risk, or safety-critical systems, tools for verification are important for ensuring correctness before testing, implementation, or deployment.  
As autonomy grows in prevalence, and systems continue to grow in size and complexity, there is a need to extend such tools to accommodate learning-enabled components in autonomous systems that resist traditional modeling.  
Stochastic reachability is an established tool for model-based verification and probabilistic safety \cite{abate2008probabilistic, summers2010verification} that has been applied to safety-critical problems in a variety of domains, including spacecraft rendezvous and docking \cite{lesser2013stochastic}, robotics and path planning \cite{malone17TRO}, and vehicle control \cite{kariotoglu2017AI}.  
Safety refers to the ability of trajectories of the system to respect known constraints on the state space, with at least a desired likelihood, despite bounded control authority. \emph{In this paper, we characterize confidence bounds on a data-driven approach for stochastic reachability, enabling rigorous assurances of safety in a model-free manner.}

Model-based stochastic reachability has received considerable attention, with methods and tools specific to Markov decision processes \cite{soudjani2015fau,cauchi2019stochy}, polynomial dynamical systems \cite{prajna2007framework, prajna2004safety}, and linear dynamical systems with log-concave disturbances \cite{VinodHSCC2017,vinod2021TAC,vinod2019sreachtools}, as well as the development of various benchmarks for comparisons \cite{abate2020arch,abate2019arch,abate2018arch}.
Additional work has been done in the model checking community to develop tools for probabilistic, bounded-time reachability problems on parametric stochastic systems which provide interval-based assurances on the safety probabilities \cite{shmarov2015probreach, shmarov2020automated}.
However, far less work has been done on data-driven stochastic reachability.  
Methods for statistical verification \cite{roohi17HSCC, zarei2020HSCC,wang19FMSD}, forward reachable set estimation \cite{thorpe2021learning, LewPavone2020}, and data-driven interval-based methods based in adaptive sampling and Gaussian processes \cite{devonport2020data} have been explored, with application to powertrains, bipedal robots, and learning-enabled cyber-physical systems.  

We focus in particular on a framework for stochastic reachability that is based in dynamic programming \cite{summers2010verification,abate2008probabilistic}.  In contrast to model checking approaches, methods based on this framework are amenable to simultaneous controller synthesis, and the underlying theory accommodates a wide range of dynamical systems that can be captured as a Markov decision process with both discrete- and continuous-valued variables.
This approach is general, in that it encompasses formulations for the terminal-hitting time problem \cite{summers2010verification, abate2008probabilistic}, as well as for the first-hitting time problem, which is far more computationally complex \cite{summers2010verification}. 
We employ a point-based formulation, in which we seek to find the likelihood of maintaining desirable time-varying state constraints over a finite time horizon, from a known initial condition.  In contrast to set-based methods \cite{gleason2017underapproximation, gleason2021lagrangian, devonport2020data}, we do not seek to compute the stochastic reachable set, the largest set of states for which the desired likelihood of safety can be maintained.



We draw in particular on tools based in reproducing kernel Hilbert spaces (RKHS), a family of machine learning techniques which enable an implicit approximation of the underlying stochastic kernel.  In our previous work \cite{thorpe2019model}, we integrated the RKHS framework with dynamic programming-based approaches for stochastic reachability, to facilitate a data-driven approach to stochastic reachability \cite{thorpe2019model,thorpe2021approximate}.  Key to the utility of this approach are bounds on the approximation error associated with the RKHS methods.  In \cite{thorpe2019model}, we described the existence of an upper bound on the asymptotic convergence rate, which can be used to infer the existence of a confidence bound on the stochastic reachability probability.  However, these bounds are difficult to compute, and hence are not practicable for ascertaining the quality of the result.  In this paper, we present \emph{computable confidence bounds} that are state- and input- dependent.  The benefit of this approach is that it enables higher accuracy bounds in regions of the state- and input- space for which more training data exists.

Kernel methods have long been used in probability and statistics, \cite{berlinet2004reproducing, parzen1961approach}, and more recently applied to
%
Markov models \cite{grunewalder2012modelling},
partially observable systems \cite{song2010hilbert, nishiyama2012hilbert}, and
policy synthesis \cite{lever2015modelling}
in the context of optimal control. 
Finite sample bounds for conditional distribution embeddings in these contexts are presented in
\cite{song2009hilbert, song2010hilbert, song2010nonparametric, grunewalder2012modelling, grunewalder2012conditional, kanagawa2014recovering},
which show that the estimated expectation
converges in probability to the true expected value with a known, asymptotic convergence rate. 
The main drawback of kernel-based approaches is the tradeoff that can occur between computability and accuracy. Higher volumes of data are required to ensure an accurate result, which inevitably increases the computational burden associated with inversion of a matrix whose size is proportional to the sample size \cite{rahimi2007random, muandet2017kernel}. Fortunately, several methods have been developed to significantly reduce the computational overhead \cite{le2013fastfood, rahimi2007random, thorpe2021approximate}.  

%
%


In this paper, we construct finite sample bounds for kernel embedding-based computation of the stochastic reachability probability, that are specific to the stochastic reachability problem.  
Our main contribution is the construction of computable state- and input-based upper and lower bounds, obtained via concentration inequalities \cite{mcdiarmid1989method} and tools from statistical learning theory \cite{vapnik1998statistical}.
Our proposed approach  
accommodates the fact that the gathered data may not be available uniformly through the state and input space.
That is, the bound is dependent upon the sample set from which the kernel embedding is inferred.  


The paper is organized as follows.
In Section \ref{section: preliminaries}, we formulate the problem and provide relevant background information. 
In Section \ref{section: finite sample bounds}, we derive the finite sample bounds.
We discuss implications of the proposed bounds in Section \ref{section: discussion}, and the problem of parameter selection.
In Section \ref{section: numerical examples}, we numerically validate the bounds on a stochastic chain of integrators, and demonstrate its utility on a nonlinear pendulum system and a nonlinear cart-pole system with black-box neural network controllers. 
Concluding remarks are provided in Section \ref{section: conclusion}.


\section{Preliminaries}
\label{section: preliminaries}

We employ the following notational conventions. 
Let $E$ be an arbitrary nonempty space.
The indicator function
$\boldsymbol{1}_{A} : E \rightarrow \lbrace 0, 1 \rbrace$
of $A \subseteq E$ is defined such that
$\boldsymbol{1}_{A}(x) = 1$ if $x \in A$, and
$\boldsymbol{1}_{A}(x) = 0$ if $x \notin A$.
%
Let $\mathcal{E}$ denote the $\sigma$-algebra on $E$.
If $E$ is a topological space \cite{ccinlar2011probability},
the $\sigma$-algebra generated by the set of all open subsets of $E$
is called the Borel $\sigma$-algebra, denoted by $\mathscr{B}(E)$.
Let $(\Omega, \mathcal{F}, \mathbb{P})$ denote a probability space,
where $\mathcal{F}$ is the $\sigma$-algebra on $\Omega$ and
$\mathbb{P} : \mathcal{F} \rightarrow [0, 1]$ is a \emph{probability measure}
on the measurable space $(\Omega, \mathcal{F})$.
A measurable function $X : \Omega \rightarrow E$ is called a \emph{random variable} taking values in $(E, \mathcal{E})$. The image of $\mathbb{P}$ under $X$, $\mathbb{P}(X^{-1}A)$, $A \in \mathcal{E}$ is called the \emph{distribution} of $X$.
Let $T$ be an arbitrary set, and for each $t \in T$, let $X_{t}$ be a random variable.
The collection of random variables $\lbrace X_{t} : t \in T \rbrace$
on $(\Omega, \mathcal{F})$ is a \emph{stochastic process}.
%


\subsection{System Model \& Stochastic Reachability Problem}

Consider a discrete-time stochastic dynamical system described by a Markov control process
as defined in \cite{summers2010verification}.

\begin{defn}[Markov Control Process]
  \label{defn: markov control process}
  A Markov control process $\mathcal{H}$ is defined as a 3-tuple,
  $\mathcal{H} = (\mathcal{X}, \mathcal{U}, Q)$,
  consisting of:
	    a Borel space $\mathcal{X} \subseteq \mathbb{R}^{n}$ representing the state space,
	    a compact Borel space $\mathcal{U} \subset \mathbb{R}^{m}$  representing the control space,
	    and $Q : \mathscr{B}(\mathcal{X}) \times \mathcal{X} \times \mathcal{U} \rightarrow [0, 1]$,
			a Borel-measurable stochastic kernel
	    which assigns a probability measure $Q(\cdot \mid x, u)$ to every $x \in \mathcal{X}$ and $u \in \mathcal{U}$ on the Borel space $(\mathcal{X}, \mathscr{B}(\mathcal{X}))$.
\end{defn}

The system evolves from an initial condition $x_{0} \in \mathcal{X}$ over a finite time horizon $t = 0, 1, \ldots, N$, $N < \infty$.
The control inputs are chosen according to a Markov control policy
$\pi = \lbrace \pi_{0}, \pi_{1}, \ldots, \pi_{N-1} \rbrace$,
which is a sequence of universally-measurable \cite{bertsekas1978stochastic} maps $\pi_{i} : \mathcal{X} \rightarrow \mathcal{U}$, $i = 0, 1, \ldots, N-1$.

For simplicity, we assume that the policy $\pi$ is stationary, meaning $\pi(x) = \pi_{0}(x) = \pi_{1}(x) = \cdots = \pi_{N-1}(x)$ for all $x \in \mathcal{X}$.
This is a simplifying assumption for the purpose of analysis, and the extension to non-stationary policies is trivial.

We consider a probabilistic reach/avoid problem known as the \emph{terminal-hitting time problem} where the objective is to determine the likelihood that the system starting at an initial condition $x_{0} \in \mathcal{X}$ will remain in a pre-defined safe set and reach a target set at the final time instant.
Let $\mathcal{K}, \mathcal{T} \in \mathscr{B}(\mathcal{X})$, denote the safe set and target set, respectively.
From \cite{summers2010verification}, the \emph{terminal-hitting time safety probability} $r_{x_{0}}^{\pi}(\mathcal{K}, \mathcal{T})$ is defined as the probability that a system $\mathcal{H}$ following a Markov policy $\pi$ from the initial condition $x_{0}$ will reach the target set $\mathcal{T}$ at time $N$ while remaining within the safe set $\mathcal{K}$ for all $t = 0, 1, \ldots, N - 1$,
\begin{equation}
	\label{eqn: terminal-hitting safety probabilities}
    r_{x_{0}}^{\pi}(\mathcal{K}, \mathcal{T}) =
    \mathbb{E}_{x_{0}}^{\pi} \Biggl[
        \Biggl(
            \prod_{i=0}^{N-1} \boldsymbol{1}_{\mathcal{K}}(x_{i})
        \Biggr)
        \boldsymbol{1}_{\mathcal{T}}(x_{N})
    \Biggr]. 
\end{equation}
For a fixed Markov policy $\pi$,
we define the value functions
$V_{t}^{\pi} : \mathcal{X} \rightarrow [0, 1]$, $t = 0, \ldots, N$,
via \emph{backward recursion}:
\begin{align}
	\label{eqn: terminal-hitting value function time N}
  V_{N}^{\pi}(x) &=
	\boldsymbol{1}_{\mathcal{T}}(x), \\
	\label{eqn: terminal-hitting value function time k}
  V_{t}^{\pi}(x) &=
	\boldsymbol{1}_{\mathcal{K}}(x)
  \int_{\mathcal{X}} V_{t+1}^{\pi}(y) Q(\mathrm{d} y \mid x, \pi(x)).
\end{align}
Then according to \cite{summers2010verification}, we have that $V_{0}^{\pi}(x_{0}) = r_{x_{0}}^{\pi}(\mathcal{K}, \mathcal{T})$
for every $x_{0} \in \mathcal{X}$.

We consider the case where the stochastic kernel $Q$ is unknown, which means the integral in \eqref{eqn: terminal-hitting value function time k} becomes intractable.

\begin{assum}
    The stochastic kernel $Q$ is unknown.
\end{assum}

Instead, we presume that a finite sample $\mathcal{S}$ of $M \in \mathbb{N}$ observations is available, taken independently and identically distributed (i.i.d.) from the stochastic kernel $Q$,
\begin{equation}
    \label{eqn: sample}
    \mathcal{S} = \lbrace (y_{i}, x_{i}, u_{i}) \rbrace_{i=1}^{M},
\end{equation}
where $x_{i}$ are taken i.i.d. from a probability measure with support $\mathcal{X}$, $u_{i} = \pi(x_{i})$, and
$y_{i} \sim Q(\cdot \mid x_{i}, u_{i})$.


\subsection{Conditional Distribution Embeddings}

Because $Q$ is assumed to be unknown, we cannot compute the safety probabilities $r_{x_{0}}^{\pi}(\mathcal{K}, \mathcal{T})$ directly. Thus, we seek to approximate the safety probabilities by numerically approximating the integral with respect to $Q(\cdot \mid x, \pi(x))$ in \eqref{eqn: terminal-hitting value function time k} using the sample $\mathcal{S}$. 
In order to approximate the integral of $V_{t+1}^{\pi}$ in \eqref{eqn: terminal-hitting value function time k}, we choose to represent the integral operator with respect to $Q(\cdot \mid x, \pi(x))$ as an element in a high-dimensional space of functions known as a reproducing kernel Hilbert space (RKHS). 

Define a positive-definite \cite[Definition~4.15]{christmann2008support} function $k : \mathcal{X} \times \mathcal{X} \to \mathbb{R}$ known as a kernel function and let $\mathscr{H}$ be a Hilbert space of functions of the form $\mathcal{X} \to \mathbb{R}$ equipped with an inner product $\langle \cdot, \cdot \rangle_{\mathscr{H}}$ and the induced norm $\lVert \cdot \rVert_{\mathscr{H}}$.


\begin{defn}[RKHS]
    \label{defn: rkhs}
    Let $E$ be an arbitrary space. 
    A Hilbert space $\mathscr{H}$ is a reproducing kernel Hilbert space (RKHS) if there exists a positive-definite kernel function $k : E \times E \to \mathbb{R}$ that satisfies the following properties~\cite{aronszajn1950theory}:
    \begin{enumerate}
        \item 
        $k(x, \cdot) \in \mathscr{H}$ for all $x \in E$, and 
        \item 
        $f(x) = \langle f, k(x, \cdot) \rangle_{\mathscr{H}}$ for all $f \in \mathscr{H}$ and $x \in E$.
    \end{enumerate}
    where the second property is known as the \emph{reproducing property}. 
\end{defn}

\begin{rem}
    Conversely, by the Moore-Aronszajn theorem \cite{aronszajn1950theory}, for any positive-definite kernel $k$, there exists a unique RKHS with $k$ as its reproducing kernel, where $\mathscr{H}$ is defined as the closure of the linear span of kernel functions, i.e. $\mathscr{H} = \overline{\mathrm{span}\lbrace k(x, \cdot) \mid x \in \mathcal{X} \rbrace}$.
    In short, this means that by defining a reproducing kernel, we obtain a corresponding RKHS.
    See, e.g., \cite{scholkopf2001learning, christmann2008support, berlinet2004reproducing} for more information on reproducing kernel Hilbert spaces.
\end{rem}

The reproducing property is central to our approach, since it allows us to evaluate any function in the RKHS as a Hilbert space inner product. By embedding the integral operator with respect to $Q(\cdot \mid x, \pi(x))$ as an element in the RKHS, we can use the reproducing property to evaluate the integral in \eqref{eqn: terminal-hitting value function time k}. 

For the measurable space $\mathcal{X}$, define the kernel function $k : \mathcal{X} \times \mathcal{X} \to \mathbb{R}$ with the associated RKHS $\mathscr{H}$. We impose a mild simplifying assumption for the purpose of analysis and assume that the kernel $k$ is bounded above by a real number $\rho < \infty$, where $\sup_{x \in \mathcal{X}} (k(x, x))^{1/2} \leq \rho < \infty$. Such a kernel is called a bounded kernel \cite[\S 4.3]{christmann2008support}.
This assumption (along with the measurability of $k$) ensures that we can represent the integral operator as an element in the RKHS.

Let $\mathscr{P}$ denote the set of probability measures on $\mathcal{X}$ which are densities of $Y \in \mathcal{X}$ conditioned on $(X, U) \in \mathcal{X} \times \mathcal{U}$ (of which the probability measures defined by the stochastic kernel $Q$ are a part). 
For any probability measure $Q(\cdot \mid x, u) \in \mathscr{P}$, if the following necessary and sufficient condition $\mathbb{E}_{Y \sim Q(\cdot \mid x, u)}[k(Y, Y)] < \infty$ is satisfied \cite{sriperumbudur2010hilbert} (which is satisfied if $k$ is bounded and measurable on $\mathcal{X}$), there exists an element $m(x, u) \in \mathscr{H}$ called the \emph{conditional distribution embedding} \cite{song2009hilbert}, defined as:
\begin{align}
    \label{eqn: kernel distribution embedding definition}
    \begin{split}
        m : \mathscr{P} & \to \mathscr{H}, \\
        Q(\cdot \mid x, u) & \mapsto m(x, u) \coloneqq \int_{\mathcal{X}} k(y, \cdot) Q(\mathrm{d} y \mid x, u).
    \end{split}
\end{align}
Using this representation, we can embed the integral with respect to $Q(\cdot \mid x, u)$ as an element in the RKHS $\mathscr{H}$.
Furthermore, if the kernel is \emph{characteristic} \cite{sriperumbudur2011universality}, then the embedding is unique, meaning the embedding captures all statistical information of the underlying distribution, and no information is lost by this representation.

As shown in \cite{thorpe2019model}, we can use conditional distribution embeddings to solve the stochastic reachability problem. By the reproducing property, we can evaluate the integral of any function $f \in \mathscr{H}$ with respect to $Q(\cdot \mid x, \pi(x))$ as a Hilbert space inner product with the embedding $m(x, \pi(x))$. 
Thus, assuming the value functions $V_{t}^{\pi}$, $t = 1, \ldots, N$, are in $\mathscr{H}$, we can evaluate the integral in \eqref{eqn: terminal-hitting value function time k} as:
\begin{align}
    & \langle V_{t+1}^{\pi}, m(x, \pi(x)) \rangle_{\mathscr{H}} \nonumber \\
    &\qquad = \biggl\langle V_{t+1}^{\pi}, \int_{\mathcal{X}} k(y, \cdot) Q(\mathrm{d} y \mid x, \pi(x)) \biggr\rangle_{\mathscr{H}} \\
    &\qquad = \int_{\mathcal{X}} \langle V_{t+1}^{\pi}, k(y, \cdot) \rangle_{\mathscr{H}} Q(\cdot \mid x, \pi(x)) \\
    &\qquad = \int_{\mathcal{X}} V_{t+1}^{\pi}(y) Q(\cdot \mid x, \pi(x)).
\end{align}
Using the embedding $m(x, \pi(x))$, we can substitute the integral expression in the stochastic reachability backward recursion \eqref{eqn: terminal-hitting value function time k} with a Hilbert space inner product in order to compute the safety probabilities $r_{x_{0}}^{\pi}(\mathcal{K}, \mathcal{T})$. This reduces the evaluation of potentially expensive integrals to a simple linear operation in Hilbert space.


\subsection{Empirically Approximating the Value Functions}

In practice, we do not have access to the true embedding $m(x, \pi(x))$ since the stochastic kernel $Q$ is unknown, meaning we cannot compute \eqref{eqn: terminal-hitting value function time k} directly. 
Thus, we compute an empirical approximation of the embedding $m(x, \pi(x))$ using a finite sample $\mathcal{S}$ as in \eqref{eqn: sample} of size $M \in \mathbb{N}$ collected i.i.d. from $Q$. 

Following \cite{micchelli2005learning}, using a sample $\mathcal{S}$, 
we can compute an empirical estimate $\hat{m}(x, u)$ of an embedding $m(x, u)$ as the solution to the following regularized least-squares problem \cite{micchelli2005learning, grunewalder2012conditional},
\begin{equation}
    \label{eqn: regularized least squares problem}
    \min_{\hat{m}} \frac{1}{M} \sum_{i=1}^{M} \lVert k(y_{i}, \cdot) - \hat{m}(x_{i}, u_{i}) \rVert_{\mathscr{H}}^{2} + \lambda \lVert \hat{m} \rVert_{\Gamma}^{2},
\end{equation}
where $\Gamma$ is a vector-valued RKHS \cite{micchelli2005learning} and $\lambda > 0$ is the regularization parameter.
By the representer theorem \cite{micchelli2005learning}, the solution to \eqref{eqn: regularized least squares problem} is unique and has the following form:
\begin{equation}
    \label{eqn: regularized least squares solution}
    \hat{m} = \sum_{i=1}^{M} \alpha_{i} k(x_{i}, \cdot) l(u_{i}, \cdot),
\end{equation}
where $\alpha \in \mathbb{R}^{M}$ is a vector of real-valued coefficients and $l : \mathcal{U} \times \mathcal{U} \to \mathbb{R}$ is a reproducing kernel function over $\mathcal{U}$.
For simplicity, we assume that the kernel $l$ is bounded and the bound $\rho$ of the kernel $k$ is also a bound for $l$.
By substituting \eqref{eqn: regularized least squares solution} into \eqref{eqn: regularized least squares problem} and taking the derivative with respect to $\alpha$, we obtain the following closed-form solution,
\begin{equation}
    \label{eqn: regularized least squares closed form solution}
    \hat{m}(x, u) = \Phi^{\top} W \Psi k(x, \cdot) l(u, \cdot),
\end{equation}
where $\Phi$ and $\Psi$ are known as \emph{feature vectors}, with elements $\Phi_{i} = k(y_{i}, \cdot)$ and $\Psi_{i} = k(x_{i}, \cdot) l(u_{i}, \cdot)$, respectively, and $W = (\Psi \Psi^{\top} + \lambda M I)^{-1}$.
For simplicity, let 
\begin{equation}
    \label{eqn: beta}
    \beta(x, u) \coloneqq W \Psi k(x, \cdot) l(u, \cdot)
\end{equation}
be a vector of real-valued coefficients that depends on the value of the conditioning variables $x$ and $u$, such that $\hat{m}(x, u) = \Phi^{\top} \beta(x, u)$.

As shown in \cite{thorpe2019model}, using an estimate $\hat{m}(x, \pi(x))$ of $m(x, \pi(x))$, we can approximate the integral of the value functions with respect to $Q(\cdot \mid x, \pi(x))$ in \eqref{eqn: terminal-hitting value function time k} as an inner product, 
\begin{align}
    \langle V_{t+1}^{\pi}, \hat{m}(x, \pi(x)) \rangle_{\mathscr{H}} 
    &= \boldsymbol{V}_{t+1}^{\pi}{}^{\top} \beta(x, \pi(x)) \nonumber \\
    &\approx \int_{\mathcal{X}} V_{t+1}^{\pi}(y) Q(\mathrm{d} y \mid x, \pi(x)),
\end{align}
where $\boldsymbol{V}_{t+1}^{\pi} = [V_{t+1}^{\pi}(y_{1}), \ldots, V_{t+1}^{\pi}(y_{M})]^{\top}$.
Thus, using the estimate $\hat{m}(x, \pi(x))$, we can recursively approximate and substitute the value functions in the stochastic reachability backward recursion in order to approximate the safety probabilities.
We summarize this procedure as Lemma \ref{lem: terminal-hitting time problem}.

\begin{lem}[Approximate Backward Recursion \cite{thorpe2019model}]
  \label{lem: terminal-hitting time problem}
	Let $\pi$ be a fixed Markov policy.
  Define the approximate value functions
  $\bar{V}_{t}^{\pi} : \mathcal{X} \rightarrow [0, 1]$, $t = 0, \ldots, N$
  by the backward recursion:
  \begin{align}
      \bar{V}_{N}^{\pi}(x) &= V_{N}^{\pi}(x), \\
      \bar{V}_{t}^{\pi}(x) &=
  	\boldsymbol{1}_{\mathcal{K}}(x)
  	\langle \bar{V}_{t+1}^{\pi}, \hat{m}_{Y|x, \pi(x)} \rangle_{\mathscr{H}}.
  \end{align}
  Then $r_{x_{0}}^{\pi}(\mathcal{K}, \mathcal{T}) \approx \bar{V}_{0}^{\pi}(x_{0})$.
\end{lem}

Note that since the estimate $\hat{m}(x, \pi(x))$ is conditioned on a particular value of $x \in \mathcal{X}$, by approximating the safety probability $r_{x_{0}}^{\pi}(\mathcal{K}, \mathcal{T})$ using Lemma \ref{lem: terminal-hitting time problem}, we obtain a point-based approximation of the safety probability at a particular value of $x_{0} \in \mathcal{X}$.
Lemma \ref{lem: terminal-hitting time problem} provides a model-free approach to approximate the stochastic reachability probability, and can easily be extended to solve related problems, including the first-hitting time problem \cite{summers2010verification} and the multiplicative and maximal cost stochastic reachability problems in \cite{abate2008probabilistic}.


\section{Finite Sample Bounds}
\label{section: finite sample bounds}

The difficulty in finding bounds on the stochastic reachability probability stems from the underlying structure of the conditional distribution embedding estimate. Unlike the embedding for a marginal distribution \cite{smola2007hilbert},
which has uniform coefficients $1/M$,
the conditional distribution embedding has non-uniform coefficients $\beta(x, u)$ \eqref{eqn: beta} which depend upon the value of the conditioning variables.
This complicates the application of existing mathematical techniques from statistical learning theory.

It is worth noting that in our case, we directly bound the variation of the estimator \emph{at a particular value of the conditioning variables $x$ and $u$}.  This means the bounds we derive provide a localized result, which can be used to assess the quality of the approximation at a particular point. This is not a significant limitation in the context of Lemma \ref{lem: terminal-hitting time problem}, since in the case of the stochastic reachability backward recursion, we seek to evaluate the safety probability at a single point.

In order to determine a bound on the quality of the approximation obtained using Lemma \ref{lem: terminal-hitting time problem}, we seek a bound on the difference between the expectation of the value function and its empirical counterpart.


\subsection{Worst-Case Difference Between the True and the Empirical Expectation}

Assume that $V_{t}^{\pi} \in \mathscr{H}$, $t = 1, \ldots, N$, and assume that for all $f \in \mathscr{H}$, $f \in [0, 1]$ and $\lVert f \rVert_{\mathscr{H}} \leq 1$.
We begin by upper bounding the deviation of the empirical expectation computed using an estimate $\hat{m}(x, u)$ from the true expectation for any value function $V_{t}^{\pi}$. 
In other words, we seek a state-based bound $B(x, u) \in \mathbb{R}$ such that for any value function $V_{t}^{\pi}$, $t = 1, \ldots, N$,
\begin{equation}
    \lvert \mathbb{E}_{Y \sim Q(\cdot \mid x, u)} [V_{t}^{\pi}(Y)] - \boldsymbol{V}_{t}^{\pi}{}^{\top} \beta(x, u) \rvert \leq B(x, u).
\end{equation}
For simplicity of notation, let $\mathbb{E} f \coloneqq \mathbb{E}_{Y \sim Q(\cdot \mid x, u)}[f(Y)]$.

We can uniformly bound the difference between the value function expectation and the empirical expectation computed using an estimate $\hat{m}(x, u)$.
Note that for any value function $V_{t}^{\pi} \in \mathscr{H}$, not necessarily attaining the supremum, 
\begin{equation}
    \label{eqn: upper bound by supremum}
    \lvert \mathbb{E} V_{t}^{\pi} - \boldsymbol{V}_{t}^{\pi}{}^{\top} \beta(x, u) \rvert \leq \sup_{\lVert h \rVert_{\mathscr{H}} \leq 1} \lvert \mathbb{E} h - \boldsymbol{h}^{\top} \beta(x, u) \rvert,
\end{equation}
where $h$ is some function in $\mathscr{H}$, $\boldsymbol{h} = [h(y_{1}), \ldots, h(y_{M})]^{\top}$, and $\beta(x, u)$ is defined as in \eqref{eqn: beta}.
This means that in the worst case, the deviation of the estimated value function expectation from the true value function expectation is less than or equal to the deviation of a function $h^{*} \in \mathscr{H}$ which satisfies the supremum. 

We then bound the right-hand side of \eqref{eqn: upper bound by supremum} using McDiarmid's inequality. In simple terms, McDiarmid's inequality states that if the empirical estimate $\boldsymbol{h}^{\top} \beta(x, u)$ computed using $\mathcal{S}$ has bounded variation when a single observation in the sample $\mathcal{S}$ is changed, then the deviation of the estimate from the true expectation is bounded by some quantity that depends on the variation bound.

\begin{defn}[Bounded Differences Condition]
	\label{defn: bounded differences condition}
	Given coefficients $c_{i} \geq 0$, $i = 1, \ldots, M$, a function $f : E^{M} \rightarrow \mathbb{R}$ satisfies the bounded differences condition if
	\begin{equation}
		\label{eqn: bounded differences condition}
		\sup_{{\substack{x_{1}, \ldots, x_{M} \\ x_{i}^{\prime} \in E}}} \vert f(x_{1}, \ldots, x_{M})
		- f(x_{1}, \ldots, x_{i}^{\prime}, \ldots, x_{M}) \vert \leq c_{i}
	\end{equation}
	for every $i = 1, \ldots, M$.
\end{defn}

\begin{lem}[McDiarmid's inequality \cite{mcdiarmid1989method}]
	\label{lem: mcdiarmids inequality}
	Let $X = \lbrace X_{1}, \ldots, X_{M} \rbrace$ be independent random variables taking values in a set $E$, and assume that the function $f : E^{M} \rightarrow \mathbb{R}$ satisfies the bounded differences condition (Definition \ref{defn: bounded differences condition}).
	Then for every $\varepsilon > 0$,
	\begin{equation}
		\label{eqn: mcdiarmids inequality}
		\Pr ( \vert f(X)
		- \mathbb{E} [f(X)] \vert \geq \varepsilon )
		\leq \exp \biggl( - \frac{2 \varepsilon^{2}}{\sum_{i=1}^{M} c_{i}^{2}} \biggr).
	\end{equation}
\end{lem}

Alternatively, Lemma \ref{lem: mcdiarmids inequality} implies that, given a small probability $\delta/2 \in (0, 1)$, then with probability $1 - \delta/2$, the deviation of the function $f$ from the expectation $\mathbb{E} [f(X)]$ is bounded by:
\begin{equation}
    f(X) - \mathbb{E} [f(X)] \leq \sqrt{ \frac{M C^{2} \log(2/\delta) }{2}}, 
\end{equation}
where $C \geq c_{i}$ for all $i = 1, \ldots, M$.

In order to determine the bound on \eqref{eqn: upper bound by supremum}, we seek 
some constant $C$ that satisfies the bounded differences condition.
However, the effect of changing an individual observation in the empirical expectation term in \eqref{eqn: upper bound by supremum} is non-trivial, since changing a single observation affects all elements of the coefficient vector $\beta(x, u)$.


%

Therefore, in order to determine $C$, we make use of a well-known result in linear algebra known as the push-through identity \cite[Fact~2.16.16]{bernstein2009matrix}.

\begin{lem}[Push-Through Identity, {\cite[Fact~2.16.16]{bernstein2009matrix}}]
    \label{lem: push through identity}
     Let $A \in \mathbb{R}^{n \times m}$ and $B \in \mathbb{R}^{m \times n}$, and assume that $AB + I$ is non-singular. Then $BA + I$ is non-singular and
    \begin{equation}
        (AB + I)^{-1} A = A(BA + I)^{-1}.
    \end{equation}
\end{lem}

We now prove the following theorem:

\begin{thm}
    \label{thm: maximum variation}
    The variation of $\boldsymbol{h}^{\top} \beta(x, u)$ by changing a single observation is at most $\rho/(\lambda M)$.
\end{thm}

\begin{pf}
    For some $h \in \mathscr{H}$, let $\boldsymbol{h}^{\top} \beta(x, u)$ be defined as in \eqref{eqn: upper bound by supremum}, such that 
    $\beta(x, u) = W \Psi k(x, \cdot) l(u, \cdot)$, where $\Psi$ is a feature vector with elements $\Psi_{i} = k(x_{i}, \cdot) l(u_{i}, \cdot)$, $W = (\Psi \Psi^{\top} + \lambda M I)^{-1}$, and the kernels $k$ and $l$ are bounded by some constant $\rho < \infty$.

    Using Lemma \ref{lem: push through identity}, we can write $W \Psi$ as:
    \begin{equation}
        \label{eqn: maximum variation identity}
        W \Psi = (\Psi \Psi^{\top} + \lambda M I)^{-1} \Psi = \Psi (\Psi^{\top} \Psi + \lambda M)^{-1},
    \end{equation}
    where we note that $(\Psi^{\top} \Psi + \lambda M)^{-1}$ is scalar.
    Using the identity in \eqref{eqn: maximum variation identity}, we can rewrite the estimate $\boldsymbol{h}^{\top} \beta(x, u)$ as:
    \begin{equation}
        \label{eqn: maximum variation proof}
        \boldsymbol{h}^{\top} W \Psi k(x, \cdot) l(u, \cdot) =
        \frac{1}{\Psi^{\top} \Psi + \lambda M} \boldsymbol{h}^{\top} \Psi k(x, \cdot) l(u, \cdot).
    \end{equation}
    Since $h \in [0, 1]$ by assumption, we have from \eqref{eqn: maximum variation proof} that the variation for changing a single observation is at most $\rho/(\lambda M)$, which proves the result.
    \qed
\end{pf}

\begin{rem}
    Note that the result provided by Theorem \ref{thm: maximum variation} is strong, and is verified by the result presented in \cite[Theorem~22]{bousquet2002stability} (with a scaling factor of $1$ instead of $1/2$, which comes from the worst-case bound).
\end{rem}

From Theorem \ref{thm: maximum variation}, we have that $C = \rho/(\lambda M)$ satisfies the bounded differences condition (Definition \ref{defn: bounded differences condition}).
Continuing from \eqref{eqn: upper bound by supremum}, using McDiarmid's inequality \eqref{eqn: mcdiarmids inequality}, we have that given $\delta/2 \in (0, 1)$, with probability $1 - \delta/2$,
\begin{align}
    & \sup_{\lVert h \rVert_{\mathscr{H}} \leq 1} \lvert \mathbb{E} h - \boldsymbol{h}^{\top} \beta(x, u) \rvert \nonumber \\
    \label{eqn: pf mcdiarmids 1}
    & \qquad \leq \mathbb{E}_{\mathcal{S}} \Biggl[ \sup_{\lVert h \rVert_{\mathscr{H}} \leq 1} \lvert \mathbb{E} h - \boldsymbol{h}^{\top} \beta(x, u) \rvert \Biggr] + \sqrt{\frac{\rho^{2} \log(2/\delta) }{2 \lambda^{2} M}}.
\end{align}
We now have an expression for the worst-case difference between the true and empirical stochastic reachability probabilities. However, because the expectation on the right-hand side of \eqref{eqn: pf mcdiarmids 1} relies upon the true expectation of $h$, it is not directly computable. 


\subsection{Removing Reliance Upon the True Expectation}

In order to enable computability of the bound in \eqref{eqn: pf mcdiarmids 1}, we bound the first term on the right-hand side of \eqref{eqn: pf mcdiarmids 1} via symmetrization and then utilize the properties of the RKHS to bound the expectation of the worst-case empirical estimate.
In particular, we make use of the reproducing property and the definition of the dual norm for Hilbert spaces, which we present here adapted from \cite[Theorem~4.3]{rudin1991functional}.

\begin{defn}[Dual Norm]
    \label{defn: dual norm}
    Let $\mathscr{H}$ be a Hilbert space. For any $f, g \in \mathscr{H}$, $\lVert f \rVert_{\mathscr{H}} = \sup_{\lVert g \rVert_{\mathscr{H}} \leq 1} \lvert \langle f, g \rangle_{\mathscr{H}} \rvert$.
\end{defn}

In addition, we rely upon a special type of random variable known as a Rademacher variable (cf. \cite{bartlett2002rademacher}), which is a uniform random variable taking values in $\lbrace -1, 1 \rbrace$.

\begin{defn}[Rademacher Variable \cite{bartlett2002rademacher}]
    A random variable $\sigma$ is called a \emph{Rademacher variable} if it is independent uniform, such that $\Pr(\sigma = 1) = \Pr(\sigma = -1) = 1/2$. 
\end{defn}

We now prove the following lemma:

\begin{lem}
    \label{lem: complexity}
    Let $h \in \mathscr{H}$ and $\boldsymbol{h}^{\top} \beta(x, u)$ be defined as in \eqref{eqn: upper bound by supremum}. Given $\delta/2 \in (0, 1)$, then with probability $1 - \delta/2$, 
    \begin{align}
        & \mathbb{E}_{\mathcal{S}} \Biggl[ \sup_{\lVert h \rVert_{\mathscr{H}} \leq 1} \lvert \mathbb{E} h - \boldsymbol{h}^{\top} \beta(x, u) \rvert \Biggr] \nonumber \\
        & \qquad \leq 
        2\sqrt{\mathrm{tr}(\boldsymbol{\beta}^{\top} \Phi \Phi^{\top} \boldsymbol{\beta})}
	    + 2\sqrt{\frac{\rho^{2} \log(2/\delta)}{2 \lambda^{2} M}}.
    \end{align}
\end{lem}

\begin{pf}
    We begin by bounding the first term on the right-hand side of \eqref{eqn: pf mcdiarmids 1} via symmetrization \cite{vapnik1998statistical}.
    Let $\tilde{\mathcal{S}}$ be a \emph{ghost sample}, that is, an independent copy of $\mathcal{S}$ that is drawn from the same sampling distribution as $\mathcal{S}$ \cite{vapnik1998statistical}.
    We replace the expectation in the first term on the right-hand side of \eqref{eqn: pf mcdiarmids 1} with a second empirical estimate computed using $\tilde{\mathcal{S}}$, to obtain:
    \begin{align}
        & \mathbb{E}_{\mathcal{S}} \Biggl[ \sup_{\lVert h \rVert_{\mathscr{H}} \leq 1} \lvert \mathbb{E} h - \boldsymbol{h}^{\top} \beta(x, u) \rvert \Biggr] \nonumber \\
        \label{eqn: pf introduce ghost sample}
        & \qquad \leq \mathbb{E}_{\mathcal{S}} \Biggl[ \sup_{\lVert h \rVert_{\mathscr{H}} \leq 1} \lvert \mathbb{E}_{\tilde{\mathcal{S}}} [ \tilde{\boldsymbol{h}}^{\top} \tilde{\beta}(x, u) - \boldsymbol{h}^{\top} \beta(x, u) ] \rvert \Biggr] \\
        \label{eqn: pf convexity of supremum}
        & \qquad \leq \mathbb{E}_{\mathcal{S} \tilde{\mathcal{S}}} \Biggl[ \sup_{\lVert h \rVert_{\mathscr{H}} \leq 1} \lvert \tilde{\boldsymbol{h}}^{\top} \tilde{\beta}(x, u) - \boldsymbol{h}^{\top} \beta(x, u) \rvert \Biggr],
    \end{align}
    where $\tilde{\beta}(x, u)$ is computed using $\tilde{\mathcal{S}}$
    as in \eqref{eqn: beta} and
    $\tilde{\boldsymbol{h}} = [h(\tilde{y}_{1}), \ldots, h(\tilde{y}_{M})]^{\top}$.
    Equation \eqref{eqn: pf introduce ghost sample} follows almost surely from the properties of conditional expectations and \eqref{eqn: pf convexity of supremum} follows by the convexity of the supremum. 
    
    We next exploit the symmetry of the empirical distributions to upper bound \eqref{eqn: pf convexity of supremum}.
    Let $\boldsymbol{\sigma}$ be Rademacher variables with
    $\boldsymbol{\sigma}_{i} \in \lbrace -1, 1 \rbrace$, and
    let $\boldsymbol{\beta} \coloneqq \mathrm{diag}(\beta(x, u))$.
    Since the distribution of the difference in empirical expectations  $\tilde{\boldsymbol{h}}{}^{\top} \tilde{\beta}(x, u) -
    \boldsymbol{h}^{\top} \beta(x, u)$ is symmetric around $0$, which follows since $\boldsymbol{f}^{\top}\beta(x, u) \in [0, 1]$ for all $f \in \mathscr{H}$,
    we see that
    $\tilde{\boldsymbol{h}}{}^{\top} \tilde{\boldsymbol{\beta}} \boldsymbol{\sigma} - \boldsymbol{h}^{\top} \boldsymbol{\beta} \boldsymbol{\sigma}$ has the same distribution.
    In effect, the Rademacher variables randomly exchange observations in $\mathcal{S}$ and $\tilde{\mathcal{S}}$ with probability $1/2$. When we take the expectation over $\boldsymbol{\sigma}$, the expectations of the empirical estimates computed using $\mathcal{S}$ and $\tilde{\mathcal{S}}$ are the same.
    Using this fact, we obtain the following:
    \begin{align}
    	\label{eqn: pf introduce rademacher variables}
    	& \mathbb{E}_{\mathcal{S} \tilde{\mathcal{S}}}
    	\Biggl[ \sup_{\lVert h \rVert_{\mathscr{H}} \leq 1}
    	\lvert \tilde{\boldsymbol{h}}{}^{\top} \tilde{\beta}(x, u) - \boldsymbol{h}^{\top} \beta(x, u) \rvert \Biggr] \nonumber \\
    	& \qquad =
    	\mathbb{E}_{\mathcal{S} \tilde{\mathcal{S}} \boldsymbol{\sigma}}
    	\Biggl[ \sup_{\lVert h \rVert_{\mathscr{H}} \leq 1} \lvert \tilde{\boldsymbol{h}}{}^{\top} \tilde{\boldsymbol{\beta}} \boldsymbol{\sigma} - \boldsymbol{h}^{\top} \boldsymbol{\beta} \boldsymbol{\sigma} \rvert \Biggr] \\
    	\label{eqn: pf symmetrization argument}
    	& \qquad \leq
    	2\mathbb{E}_{\mathcal{S} \boldsymbol{\sigma}} \Biggl[ \sup_{\lVert h \rVert_{\mathscr{H}} \leq 1} \lvert \boldsymbol{h}^{\top} \boldsymbol{\beta} \boldsymbol{\sigma} \rvert \Biggr].
    \end{align}
    Then by applying the reproducing property (Definition \ref{defn: rkhs})
    and the definition of the dual norm for Hilbert spaces (Definition \ref{defn: dual norm}), we remove the dependence on $h$:
    \begin{align}
    	& 2\mathbb{E}_{\mathcal{S} \boldsymbol{\sigma}}
    	\Biggl[
    	\sup_{\lVert h \rVert_{\mathscr{H}} \leq 1}
    	\vert
    	\boldsymbol{h}^{\top} \boldsymbol{\beta} \boldsymbol{\sigma}
    	\vert
    	\Biggr] \nonumber \\
    	& \qquad =
    	2\mathbb{E}_{\mathcal{S} \boldsymbol{\sigma}}
    	\Biggl[
    	\sup_{\lVert h \rVert_{\mathscr{H}} \leq 1}
    	\vert
    	\langle h, \Phi^{\top} \boldsymbol{\beta} \boldsymbol{\sigma} \rangle_{\mathscr{H}}
    	\vert
    	\Biggr] \\
    	& \qquad =
    	2\mathbb{E}_{\mathcal{S} \boldsymbol{\sigma}}
    	[
    	\lVert
    	\Phi^{\top} \boldsymbol{\beta} \boldsymbol{\sigma}
    	\rVert_{\mathscr{H}}
    	].
    \end{align}
    We then utilize the definition of the Hilbert space norm
    and the concavity of the square root and note that the expectation of the Rademacher variables $\mathbb{E}_{\boldsymbol{\sigma}}[\boldsymbol{\sigma}_{i}\boldsymbol{\sigma}_{j}]$ vanishes except when $i = j$,
    \begin{align}
    	& 2\mathbb{E}_{\mathcal{S} \boldsymbol{\sigma}}
    	[
    	\lVert
    	\Phi^{\top} \boldsymbol{\beta} \boldsymbol{\sigma}
    	\rVert_{\mathscr{H}}
    	] \nonumber \\
    	& \qquad \leq
    	2 \mathbb{E}_{\mathcal{S}}
    	\Biggl[
    	\Bigl(
    	\mathbb{E}_{\boldsymbol{\sigma}} \Bigl[
    	\boldsymbol{\sigma}^{\top} \boldsymbol{\beta}^{\top} \Phi
    	\Phi^{\top} \boldsymbol{\beta} \boldsymbol{\sigma}
    	\Bigr]
    	\Bigr)^{1/2}
    	\Biggr] \\
    	\label{eqn: pf trace}
    	& \qquad =
    	2 \mathbb{E}_{\mathcal{S}}
    	\Biggl[
    	\Bigl(
    	\mathrm{tr}\Bigl(
    	\boldsymbol{\beta}^{\top} \Phi
    	\Phi^{\top} \boldsymbol{\beta}
    	\Bigr)
    	\Bigr)^{1/2}
    	\Biggr].
    \end{align}
    By bounding the expectation in \eqref{eqn: pf trace} by McDiarmid's inequality again, we obtain:
    \begin{align}
    	\label{eqn: pf mcdiarmids 2}
    	& 2\mathbb{E}_{\mathcal{S}}
    	\Biggl[
    	\Bigl(
    	\mathrm{tr}\Bigl(
    	\boldsymbol{\beta}^{\top} \Phi
    	\Phi^{\top} \boldsymbol{\beta}
    	\Bigr)
    	\Bigr)^{1/2}
    	\Biggr] \nonumber \\
    	& \quad \leq
    	2\sqrt{\mathrm{tr}(\boldsymbol{\beta}^{\top} \Phi \Phi^{\top} \boldsymbol{\beta})}
    	+ 2\sqrt{\frac{\rho^{2} \log(2/\delta)}{2 \lambda^{2} M}},
    \end{align}
    which proves the result.
    \qed
\end{pf}

Continuing from \eqref{eqn: pf mcdiarmids 1}, using McDiarmid's inequality with $C = \rho/(\lambda M)$ via Theorem \ref{thm: maximum variation} and Lemma \ref{lem: complexity}, we have that given $\delta/2 \in (0, 1)$, with probability $1 - \delta/2$, 
\begin{align}
    \label{eqn: uniform upper bound}
	&
	\sup_{\lVert h \rVert_{\mathscr{H}} \leq 1}
	\vert \mathbb{E} h -
	\boldsymbol{h}^{\top}\beta(x, u) \vert \nonumber \\
	& \qquad \leq
	2\sqrt{\mathrm{tr}(\boldsymbol{\beta}^{\top} \Phi \Phi^{\top} \boldsymbol{\beta})} 
	+ 3\sqrt{\frac{ \rho^{2} \log(2/\delta)}{2 \lambda^{2} M}}.
\end{align}
Thus, we have a computable bound for the worst-case difference between the true and empirical expectation of a function $h \in \mathscr{H}$.

\begin{rem}
    We note that even in the worst case, the first term on the right-hand side of \eqref{eqn: uniform upper bound} is bounded, since
    $\lVert \Phi^{\top} \beta(x, u) \rVert_{\mathscr{H}} \leq 1$ by assumption.
    However, in practice, with appropriate kernel selection and choice of regularization parameter $\lambda$, this term will often be significantly less than $1$.
\end{rem}

\subsection{Finite Sample Bound on the Safety Probability}

We now can state the main result, which we present as Theorem \ref{thm: main result}.
\begin{thm}
    \label{thm: main result}
    For any value function $V_{t}^{\pi}$, given $\delta/2 \in (0, 1)$, with probability $1-\delta/2$,
    the difference between the true and empirical expectation of the value functions is bounded by:
    \begin{align}
        \label{eqn: finite sample bound result}
        & \lvert \mathbb{E} V_{t}^{\pi} - \boldsymbol{V}_{t}^{\pi}{}^{\top} \beta(x, u) \rvert \nonumber \\
        & \qquad \leq
        2\sqrt{\mathrm{tr}(\boldsymbol{\beta}^{\top} \Phi \Phi^{\top} \boldsymbol{\beta})}
        + 3\sqrt{\frac{\rho^{2} \log(2/\delta)}{2 \lambda^{2} M}}.
    \end{align}
\end{thm}

\begin{pf}
    The proof follows from the arguments presented for \eqref{eqn: uniform upper bound}. 
    In \eqref{eqn: upper bound by supremum},
    we uniformly bound the difference between the true and empirical expectation of the value functions by the worst-case function $h \in \mathscr{H}$.
    We then use McDiarmid's inequality (Lemma \ref{lem: mcdiarmids inequality}) with the bound $C = \rho/(\lambda M)$ satisfying the bounded differences condition (Theorem \ref{thm: maximum variation}) to obtain \eqref{eqn: pf mcdiarmids 1}. 
    Using a ghost sample, the symmetrization argument, and the definition of the dual norm for Hilbert spaces, we then bound the first term on the right-hand side of \eqref{eqn: pf mcdiarmids 1} in order to remove the dependence on $h$,
    \begin{align}
        \label{eqn: main thm step 3}
        &
    	\mathbb{E}_{\mathcal{S}}
    	\Biggl[
    	\sup_{\lVert h \rVert_{\mathscr{H}} \leq 1}
    	\vert
    	\mathbb{E} h -
    	\boldsymbol{h}^{\top} \beta(x, u)
    	\vert
    	\Biggr] \nonumber \\
    	& \quad \leq
        2\sqrt{\mathrm{tr}(\boldsymbol{\beta}^{\top} \Phi \Phi^{\top} \boldsymbol{\beta})}
    	+ 2\sqrt{\frac{\rho^{2} \log(2/\delta)}{2 \lambda^{2} M}}.
    \end{align}
    We then substitute the bound in \eqref{eqn: main thm step 3} into \eqref{eqn: pf mcdiarmids 1} and obtain the result,
    \begin{align}
        \label{eqn: main thm step 4}
        & \vert \mathbb{E} V_{t}^{\pi} - \boldsymbol{V}_{t}^{\pi}{}^{\top} \beta(x, u) \vert \nonumber \\
    	& \quad \leq
    	2\sqrt{\mathrm{tr}(\boldsymbol{\beta}^{\top} \Phi \Phi^{\top} \boldsymbol{\beta})}
    	+ 3\sqrt{\frac{\rho^{2} \log(2/\delta)}{2 \lambda^{2} M}},
    \end{align}
    which proves \eqref{eqn: finite sample bound result}. \qed
\end{pf}
Let $B(x, u)$ be the bound on the difference between the expected value of the value functions and its empirical counterpart in \eqref{eqn: finite sample bound result}, given by
\begin{equation}
    \label{eqn: bound B}
    B(x, u) = 2\sqrt{\mathrm{tr}(\boldsymbol{\beta}^{\top} \Phi \Phi^{\top} \boldsymbol{\beta})}
    	+ 3\sqrt{\frac{\rho^{2} \log(2/\delta)}{2 \lambda^{2} M}}.
\end{equation}
This means
that for any value function $V_{t}^{\pi}$, given $\delta/2 \in (0, 1)$, with probability $1-\delta/2$, that the absolute difference between the actual expectation and the empirical expectation computed using $\hat{m}(x, \pi(x))$ is bounded by 
\begin{equation}
    -B(x, u) \leq \mathbb{E} V_{t}^{\pi} - \boldsymbol{V}_{t}^{\pi}{}^{\top} \beta(x, u) \leq B(x, u).
\end{equation}

Thus, by applying this bound to the expectations in the backward recursion, we obtain an overall bound on the approximation of the safety probabilities obtained using Lemma \ref{lem: terminal-hitting time problem}. 
Furthermore, the bound in \eqref{eqn: bound B} depends on the value of the conditioning variables $x$ and $u$, which means the bound can serve as an indication of the quality of the approximation at a particular point. 

Note that the bound in \eqref{eqn: bound B} applies to the error of a value function at a single time step. This means that the total error on the approximation of the safety probabilities increases linearly with the number of time steps used in the backward recursion in Lemma \ref{lem: terminal-hitting time problem}. 
Following \cite[Corollary~3]{thorpe2019model}, if the error of each value function is bounded by $B(x, u)$ in \eqref{eqn: bound B}, we obtain an overall bound on the safety probabilities $V_{0}^{\pi}$ of $N B(x, u)$, where $N$ is the time horizon.
We summarize this result in the following lemma:

\begin{lem}
    Given a bound $B(x, u)$, the error in the safety probabilities obtained via Lemma \ref{lem: terminal-hitting time problem} is given by
    \begin{equation}
        \lvert \mathbb{E} V_{0}^{\pi} - \boldsymbol{V}_{0}^{\pi}{}^{\top} \beta(x, u) \rvert \leq N B(x, u),
    \end{equation}
    where $N$ is the time horizon.
\end{lem}

\begin{pf}
    The proof is by induction. If the error of a single value function is $B(x, u)$ as in \eqref{eqn: bound B}, and the safety probabilities are computed via recursive substitution as in Lemma \ref{lem: terminal-hitting time problem}, then at time $t = 0$, the maximum error of $V_{0}^{\pi}$ is $N B(x, u)$. 
\end{pf}


\section{Kernel and Parameter Selection}
\label{section: discussion}

The quality of the approximation of the stochastic reachability probability is governed by a number of factors, including the choice of kernel function, parameters associated with the kernel function, the regularization parameter of the least-squares problem \eqref{eqn: regularized least squares problem}, and the sample, $\mathcal{S}$.
In a realistic setting, we typically do not have explicit control over the sample $\mathcal{S}$ or the number of observations in the sample. Thus, the choice of kernel function and the model parameters plays an important role in the quality of the kernel-based approximation. 



\subsection{Kernel Selection}

The performance of kernel-based learning algorithms is closely tied to the choice of kernel function and the structure of the RKHS. In essence, the Hilbert space needs to be rich enough to model the set of probability measures underlying the observed data without overfitting.
Thus, we frame the problem of kernel selection as a problem of limiting the complexity of the RKHS \cite{vapnik1998statistical, bartlett2002rademacher}, where \emph{complexity} in statistical learning literature refers to the ability of a function class to fit random noise.
Intuitively, by choosing a function class that lowers the complexity term, we 
reduce the possibility that our function class $\mathscr{H}$ will overfit the observed data.

We propose a complexity term that is closely related to the Rademacher complexity \cite{bartlett2002rademacher} from statistical learning theory, but instead based on the bound presented in Theorem \ref{thm: main result} that accommodates non-uniform coefficients $\beta(x, u)$.
Define the random variable
\begin{equation}
	\hat{\mathscr{C}}(\mathscr{H})
	= \mathbb{E}_{\boldsymbol{\sigma}}
	\Biggl[
	\sup_{\lVert h \rVert_{\mathscr{H}} \leq 1}
	\vert \boldsymbol{h}^{\top} \boldsymbol{\beta} \boldsymbol{\sigma} \vert
	\Biggr].
\end{equation}
Then the conditional complexity of $\mathscr{H}$ is defined as $\mathscr{C}(\mathscr{H}) = \mathbb{E}_{\mathcal{S}}[\hat{\mathscr{C}}(\mathscr{H})]$. 
Using the bound in \eqref{eqn: pf trace}, we have
\begin{equation}
    \mathscr{C}(\mathscr{H}) \leq
    \mathbb{E}_{\mathcal{S}}
	\Biggl[
	\Bigl(
	\mathrm{tr}\Bigl(
	\boldsymbol{\beta}^{\top} \Phi
	\Phi^{\top} \boldsymbol{\beta}
	\Bigr)
	\Bigr)^{1/2}
	\Biggr].
\end{equation}
This choice of complexity term is equivalent to the first term on the right-hand side of \eqref{eqn: bound B}.
Since the conditional complexity appears in the bounds presented in Theorem \ref{thm: main result}, minimizing the complexity term also minimizes the finite sample bounds on the difference in expectations.
Thus, we can choose the kernel functions which minimize the complexity term for all $(x, u)$, effectively minimizing the finite sample bound on the difference in expectations of the value functions. 

The kernel should also be chosen to satisfy
universality \cite{sriperumbudur2010hilbert, sriperumbudur2011universality} (resp. characteristic) and boundedness properties.
One common choice of kernel that satisfies these properties is the Gaussian RBF kernel $k(x, x') = \exp(-\lVert x - x' \rVert_{2}^{2}/2\sigma^{2})$, where $\sigma > 0$.
Universal \cite{smola2007hilbert, sriperumbudur2010hilbert} kernels are so-named because they satisfy a universal approximation property and are able to learn any real-valued function arbitrarily well. 
Because the mapping from the set of all probability measures $\mathscr{P}$ into the RKHS \eqref{eqn: kernel distribution embedding definition} is injective for universal (resp. characteristic, see \cite{sriperumbudur2011universality}) kernels, 
this means there is a unique element in the RKHS $\mathscr{H}$
for any $\mathbb{P}, \mathbb{Q} \in \mathscr{P}$, such that
$\lVert m_{\mathbb{P}} - m_{\mathbb{Q}} \rVert_{\mathscr{H}} = 0$ if and only if $\mathbb{P} = \mathbb{Q}$. 
This ensures that the conditional distribution embedding admits a unique solution \cite{sriperumbudur2010hilbert}, and that we can distinguish between distributions in Hilbert space. 
By choosing a bounded kernel function, we ensure that $\rho < \infty$, and we can achieve tighter bounds by selecting a kernel function with small $\rho$.
The Gaussian kernel function, for example, has $\rho = 1$.

Typically, a parameterized kernel which is known to satisfy these properties is chosen, and then kernel parameters are tuned via cross-validation techniques. 
However, nascent work has posed 
kernel synthesis for a given sample $\mathcal S$ as an optimization problem.  This approach has been demonstrated via convex optimization \cite{kim2008learning} and semi-definite programming \cite{lanckriet2004learning} for marginal distributions with scalar-valued regression.
While this idea is promising, the connection between conditional distribution embeddings and the underlying regression problem has only recently begun to be explored \cite{grunewalder2012conditional}.  Further, the extension from scalar-valued regression to vector-valued regression \cite{micchelli2005learning} that is required by the objective function in \eqref{eqn: regularized least squares problem} is not straightforward.    
Hence additional work will be needed to evaluate the feasibility of kernel synthesis methods for this problem.

\subsection{Parameter Selection}

Tighter bounds may be possible by identifying a strict upper bound on the elements of $\beta(x, u)$, which in turn are influenced by the parameters $\rho$, the upper bound on the kernel function, and $\lambda$, the least-squares regularization coefficient in \eqref{eqn: regularized least squares problem}.

The value of $\rho$ is determined primarily by the choice of kernel and its parameters, and can be tuned by minimizing the complexity term $\mathscr{C}(\mathscr{H})$.  
The value of $\lambda$ affects the convergence rate associated with the uniform finite sample bound.
The convergence guarantees in \cite{song2009hilbert} and \cite{grunewalder2012modelling} typically depend upon $\lambda$ going to zero as the number of observations increases. See \cite{caponnetto2007optimal, de2005model} for a discussion of optimal values of $\lambda$.


\subsection{Scalability}

One significant advantage of using the kernel-based approach in Lemma \ref{lem: terminal-hitting time problem} is that approximating the value function expectation in \eqref{eqn: terminal-hitting value function time k} using conditional distribution embeddings does not scale exponentially as the system dimensionality increases (also known as the curse of dimensionality). This is primarily due to the fact that the system dimensionality (as well as the input dimension) only has an effect on the evaluation of the kernel functions $k$ and $l$. 
Instead, the complexity of computing the conditional distribution embedding estimate is generally $\mathcal{O}(M^{3})$, which is primarily driven by the matrix inverse $W$ in \eqref{eqn: beta}, and scales solely with the sample size $M$ used to construct the embedding estimate. 
This is demonstrated empirically in \cite{thorpe2019model}, which shows that the computational complexity increases roughly linearly as the system dimensionality is increased.
However, in order to adequately characterize the stochastic kernel of a high-dimensional state space, a large sample size may be required.
In effect, this means that the quality of the approximation is in large part governed by the sample used to characterize the region of interest--which further motivates our approach to state-based finite sample bounds.

As the time horizon $N$ increases, the backward recursion in Lemma \ref{lem: terminal-hitting time problem} shows that if the policy $\pi$ is time-invariant, meaning $\pi_{0} = \pi_{1} = \cdots = \pi_{N-1}$, we do not need to recompute the embedding estimate $\hat{m}(x, \pi(x))$ at every time step. This means we only need to compute the estimate \emph{once} for a given sample $\mathcal{S}$, and indicates that the complexity of Lemma \ref{lem: terminal-hitting time problem} increases linearly as the time horizon increases. In the case of time-varying policies, this of course means that we must recompute the embedding at every time step since the closed-loop dynamics of the system vary with time in accordance with the policy.


\section{Examples}
\label{section: numerical examples}

\begin{figure*}
    \centering
    \includegraphics[keepaspectratio]{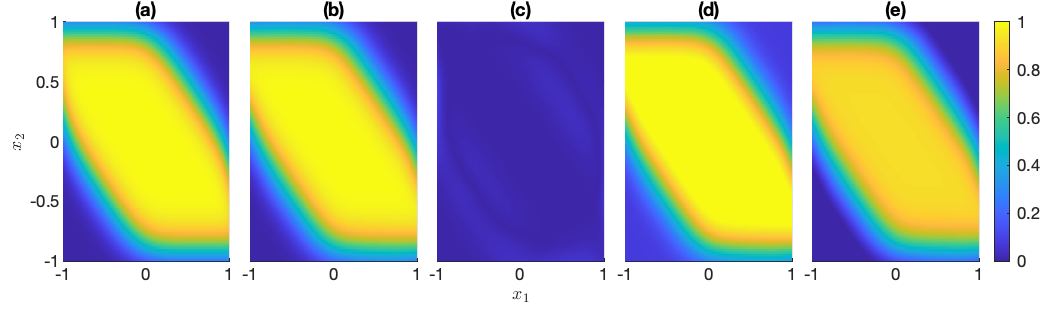}
    \caption{
    (a) Safety probabilities at $k = 0$, $N = 5$, for a double integrator computed using dynamic programming. 
    (b) Safety probabilities at $k = 0$ for a double integrator computed using Lemma \ref{lem: terminal-hitting time problem} for $N = 5$.
    (c) Absolute error $\vert V_{0}^{\pi}(x) - \bar{V}_{0}^{\pi}(x) \vert$ between the dynamic programming solution and Lemma \ref{lem: terminal-hitting time problem}.
    (d, e) Upper and lower finite sample bounds, respectively, of the safety probabilities of a double integrator system computed using Theorem \ref{thm: main result} with $\delta = 0.1$, where
    (d) is $\bar{V}_{0}^{\pi}(x) + B(x, \pi(x))$ and 
    (e) is $\bar{V}_{0}^{\pi}(x) - B(x, \pi(x))$, where $B(x, \pi(x))$ is computed as in \eqref{eqn: bound B}.
    }
    \label{fig: double integrator}
\end{figure*}

We implemented Lemma \ref{lem: terminal-hitting time problem} on a stochastic chain of integrators for the purpose of validation, and on a nonlinear pendulum system \cite{dutta2019reachability} and a closed-loop nonlinear cart-pole benchmark system \cite{manzanas2019arch} with black-box neural network controllers to demonstrate the capabilities of the proposed approach.
For each problem, we generated a sample 
$\mathcal{S}$
of observations via simulation, and then assumed no knowledge of the system dynamics or the structure of the disturbance for the purposes of computing the stochastic reachability probability.
We then computed finite sample bounds via Theorem \ref{thm: main result}.
For all problems, we chose a Gaussian kernel $k(x, x') = \exp(- \lVert x - x' \rVert_{2}^{2}/2 \sigma^{2})$ with $\sigma = 0.1$, and chose $\lambda$ using the optimal rate computed in \cite{caponnetto2007optimal}.

All computations were done in Matlab,
and code to reproduce the analysis and figures is available at:	
\texttt{github.com/unm-hscl/ajthor-Automatica2020a}.


\subsection{Stochastic Chain of Integrators}

We first consider a $2$-D stochastic chain of integrators
\cite{vinod2017scalable}, in which the input appears at the second derivative
and each element of the state vector is the discretized integral of the element that follows it. The dynamics with sampling time $T = 0.25$ are given by:
\begin{align}
	\label{eqn: double integrator dynamics}
  \boldsymbol{x}_{k+1} =
  \begin{bmatrix}
    1 & T \\
    0 & 1
  \end{bmatrix}
  \boldsymbol{x}_{k} +
  \begin{bmatrix}
    \frac{T^{2}}{2} \\
    T
  \end{bmatrix}
  u_{k} +
  \boldsymbol{w}_{k}
\end{align}
where $\boldsymbol{w}_{k}$ is an i.i.d. disturbance
defined on the probability space
$(\mathcal{W}, \mathscr{B}(\mathcal{W}), \Pr_{\boldsymbol{w}})$.
We presume a Gaussian disturbance $\boldsymbol{w}_{k} \sim \mathcal{N}(0, \Sigma)$ with $\Sigma = 0.01 I$, a control policy $\pi(x) = 0$, and target and safe sets $\mathcal{T} = [-1, 1]^{2}$ and
$\mathcal{K} = [-1, 1]^{2}$.

\begin{figure}
  \includegraphics[width=\columnwidth, height=150pt]{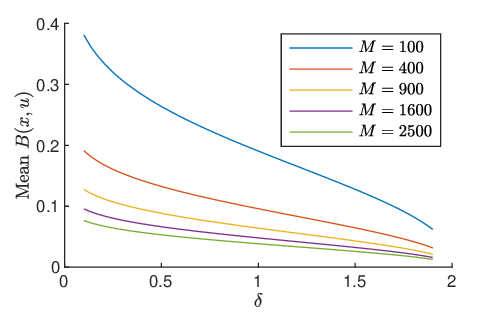}
  \caption{Figure showing the mean of the finite sample bounds $B(x, \pi(x))$ in the region $[-1, 1] \times [-1, 1]$ for a double integrator system as a function of $\delta$.}
  \label{fig: double integrator mean finite sample bounds delta}
\end{figure}

\begin{figure*}
    \centering
    \includegraphics{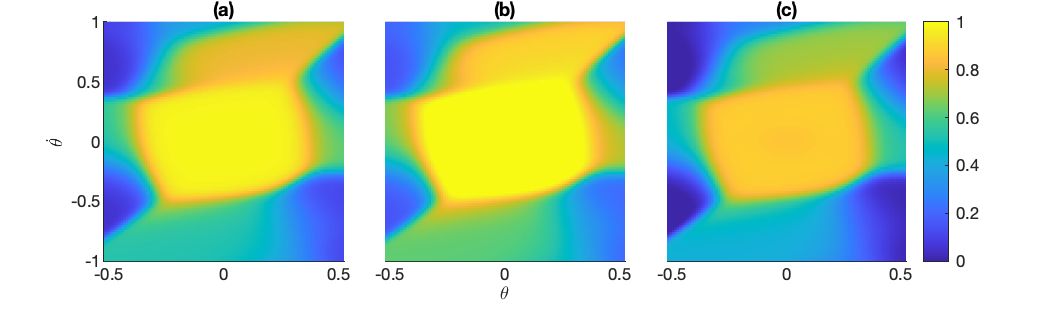}
    \caption{
    (a) Safety probabilities for a linearized cart-pole system computed using Lemma \ref{lem: terminal-hitting time problem} for $N = 10$.
    (b) Upper bound on the safety probabilities computed using Theorem \ref{thm: main result}.
    (c) Lower bound on the safety probabilities computed using Theorem \ref{thm: main result}.
    }
    \label{fig: cart pole}
\end{figure*}

We consider a sample $\mathcal{S}$ of $M = 2500$ observations drawn i.i.d. from $Q$, a representation of \eqref{eqn: double integrator dynamics} as a Markov control process  (Definition \ref{defn: markov control process}).  
The initial conditions $x \in \mathcal{X}$ in the sample were chosen uniformly in the interval $[-1.1, 1.1] \times [-1.1, 1.1]$ in order to ensure that a subset of the initial conditions violates the safety constraints, $\mathcal{K}$ and $\mathcal{T}$. We do this to ensure the ``learned'' model does not map all initial conditions to a safe set. Otherwise, in the regression, the value function estimate maps all values to $1$. The resulting state $y \in \mathcal{X}$ is drawn from $Q(\cdot \mid x, \pi(x))$ using the dynamics in \eqref{eqn: double integrator dynamics}. 
We then presumed no knowledge of the system dynamics or $Q$ and computed the estimate $\hat{m}(x, \pi(x))$ according to \eqref{eqn: regularized least squares closed form solution}, with $\beta(x, \pi(x))$ computed as in \eqref{eqn: beta}.
Using $\hat{m}(x, \pi(x))$, we then computed the stochastic reachability probability using Lemma \ref{lem: terminal-hitting time problem} for a time horizon of $N = 5$. 
We compare the stochastic reachability probability computed according to 
Lemma \ref{lem: terminal-hitting time problem}
against the solution via dynamic programming, that presumes the stochastic kernel $Q$ is known. 
The absolute error $\vert V_{0}^{\pi}(x) - \bar{V}_{0}^{\pi}(x) \vert$ between the results obtained from Lemma \ref{lem: terminal-hitting time problem} and the dynamic programming solution is shown in
Figure \ref{fig: double integrator}(c).
As expected, the stochastic reachability probabilities computed using Lemma \ref{lem: terminal-hitting time problem}, show low absolute error as compared with the dynamic programming solution, with a maximum absolute error of $0.1158$ and a mean absolute error of $0.0122$. 


We then evaluated the finite sample bounds of the approximation using Theorem \ref{thm: main result}. We computed the bound $B(x, \pi(x))$ as in \eqref{eqn: bound B} with $\delta = 0.1$ to obtain bounds on the safety probabilities. 
The bounds obtained from Theorem \ref{thm: main result} are probabilistic upper and lower bounds, meaning that given any $\delta/2 \in [0, 1]$, with probability $1 - \delta/2$, the approximation from Lemma \ref{lem: terminal-hitting time problem} is bounded by Theorem \ref{thm: main result}, which validates the result.
%
Figure \ref{fig: double integrator}(d) shows the upper bound on the safety probabilities, while Figure \ref{fig: double integrator}(e) shows the lower bound.


We can see that the difference between the upper and lower bounds is small, which indicates that the quality of the approximation obtained via Lemma \ref{lem: terminal-hitting time problem} is close to the true solution with high probability. 
Further, the absolute difference is larger than the error values in Figure \ref{fig: double integrator}(c), as expected, meaning that the computed bound is a reasonable probabilistic upper bound.

We then consider the effect of varying the parameters $M$ and $\delta$ in \eqref{eqn: bound B} on the finite sample bound computed using Theorem \ref{thm: main result}. 
In order to demonstrate the effect of the parameter $M$, we drew $5$ new samples from the stochastic kernel $Q$, where the number of observations in each sample was chosen to be of size $M \in [100, 2500]$.
An estimate was then computed for each sample and the finite sample bounds were computed for each estimate. 
For each sample of length $M \in [100, 2500]$, we computed the mean of the finite sample bounds in the region $[-1, 1] \times [-1, 1]$ for different values of $\delta \in [0.1, 1.9]$.
Figure \ref{fig: double integrator mean finite sample bounds delta} shows the mean of $B(x, \pi(x))$ \eqref{eqn: bound B}, the finite sample bounds computed using Theorem \ref{thm: main result} for the five samples of length $M \in [100, 2500]$ as a function of $\delta$.

As expected, we can see that as the size of the samples increases, we obtain tighter probabilistic bounds $B(x, \pi(x))$ in the region $[-1, 1] \times [-1, 1]$ via Theorem \ref{thm: main result}. 
Effectively, this means that as the number of observations from the stochastic kernel increases, we obtain a better estimate of the conditional distribution embedding $m(x, \pi(x))$, and thus a better estimate of the safety probabilities via Lemma \ref{lem: terminal-hitting time problem}. 
Also as expected, as the violation threshold decreases (i.e., $\delta$ increases), the mean of the probabilistic bound decreases. 
Figure \ref{fig: double integrator mean finite sample bounds delta} also shows that for low values of delta, we obtain higher values of the probabilistic bound $B(x, \pi(x))$. This corresponds to a high desired confidence.
Further, we can see in Figure \ref{fig: double integrator mean finite sample bounds delta} that the finite sample bounds do not improve appreciably as the sample size increases beyond a certain point. 



\subsection{Linearized Cart-Pole System}

We then considered a benchmark cart-pole system \cite{manzanas2019arch} with a black-box neural network feedback controller.
The dynamics for the linearized cart-pole system \cite{manzanas2019arch} are given by:
\begin{align}
\begin{split}
\label{eqn: linearized cart pole dynamics}
	\ddot{x} &= 0.0043 \dot{\theta} - 2.75 \theta + 1.94 u - 10.95 \dot{x} \\
	\ddot{\theta} &= 28.58 \theta - 0.044 \dot{\theta} - 4.44 u + 24.92 \dot{x}
\end{split}
\end{align}
with state $x = [x, \dot{x}, \theta, \dot{\theta}]^{\top} \in \mathbb{R}^{4}$ and control input $u \in \mathbb{R}$. The dynamics are then discretized in time with sampling time $T_{s} = 0.2$ s.
We add an additional Gaussian disturbance $\boldsymbol{w}_{k} \sim \mathcal{N}(0, \Sigma)$ with $\Sigma = 0.01 I$ to the dynamical state equations, which can simulate dynamical uncertainty or minor system perturbations.
The control input is computed via a neural network controller \cite{manzanas2019arch}, which takes the current state and outputs a real number $u \in \mathbb{R}$, which can be interpreted as the input torque. 

The benchmark is defined \cite{manzanas2019arch} such that the neural network controller must 
keep the lateral position of the cart $x$ within $[-0.7, 0.7]$, 
maintain a low cart velocity $\dot{x} \in [-1, 1]$, 
and keep the pendulum angle $\theta$ within $[-\pi/6, \pi/6]$ 
while the angular velocity $\dot{\theta}$ is unconstrained. 
We define the safe set $\mathcal{K}$ according to the above constraints, and define the target set $\mathcal{T}$ such that the pendulum angle $\theta$ must be within $[-0.05, 0.05]$.
\begin{align}
    \mathcal{K} &= \lbrace x \in \mathbb{R}^{4} \mid \vert x_{1} \vert \leq 0.7, \vert x_{2} \vert \leq 1, \vert x_{3} \vert \leq \pi/6 \rbrace  \\
    \mathcal{T} &= \lbrace x \in \mathbb{R}^{4} \mid \vert x_{3} \vert \leq 0.05 \rbrace
\end{align}
We simulated $10$ trajectories from initial conditions taken uniformly from the ranges specified above, and extracted a sample $\mathcal{S}$ of $M = 12{,}234$ observations taken i.i.d. from the stochastic kernel $Q$, a representation of the dynamics \eqref{eqn: linearized cart pole dynamics} as a Markov control process (Definition \ref{defn: markov control process}). 

\begin{figure}
  \includegraphics[keepaspectratio, width=\columnwidth]{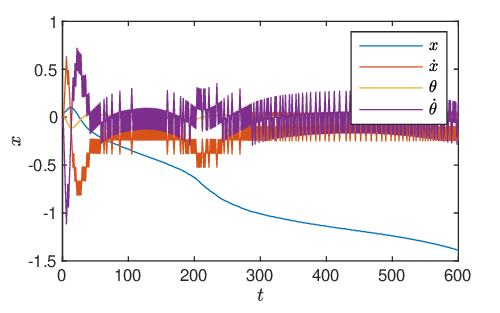}
  \caption{Sample realization of the nonlinear cart-pole system \eqref{eqn: nonlinear cartpole system} over 600 time steps.}
  \label{fig: cartpole nonlinear trajectory}
\end{figure}

We then computed the safety probabilities for the system over a time horizon $N = 10$ using Lemma \ref{lem: terminal-hitting time problem} to demonstrate the feasibility of the approach.
The safety probabilities computed using Lemma \ref{lem: terminal-hitting time problem} at $k = 0$ for $N = 10$ are shown in Figure \ref{fig: cart pole}(a). 
We can see that the closed-loop system has a high probability of stabilizing the pendulum from an initial condition within the range $\theta \in [-\pi/6, \pi/6]$ and the results also show an underlying symmetry around $\theta = 0$, as expected.
We then computed the finite sample bounds on the approximation using Theorem \ref{thm: main result} with $\delta = 0.1$ to obtain probabilistic upper and lower bounds on the safety probabilities. 
Figure \ref{fig: cart pole}(b) shows the probabilistic upper bound on the safety probabilities, while Figure \ref{fig: cart pole}(c) shows the lower bound. 
As expected, because we use a high number of observations $M = 12{,}234$, the bounds computed in Theorem \ref{thm: main result} 
show that with high probability, the solution is close to the true solution.

This means that using the proposed data-driven approach, we can utilize stochastic reachability to analyze the safety properties of a dynamical system with a black-box neural network controller. Similarly, we can expose the underlying structure of the closed-loop system to reveal useful knowledge of the system properties, such as symmetry. 



\subsection{Nonlinear Cart-Pole System}

\begin{figure}
  \includegraphics[keepaspectratio, width=\columnwidth]{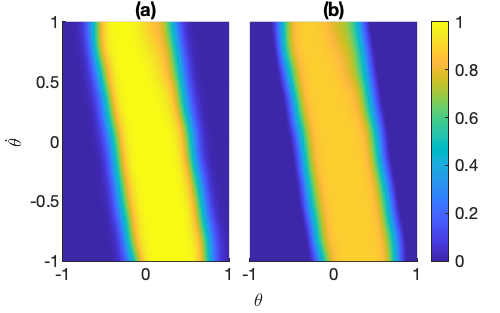}
  \caption{
  (a) Safety probabilities of the nonlinear cart-pole system computed using Lemma \ref{lem: terminal-hitting time problem} for a time horizon of $N = 600$. (b) Probabilistic lower bound on the safety probabilities computed using Theorem \ref{thm: main result} with $\delta = 0.1$. 
  }
  \label{fig: cartpole nonlinear}
\end{figure}

We then analyzed a nonlinear cart-pole system with a neural network controller \cite{manzanas2019arch}, with dynamics given by:
\begin{align}
\begin{split}
\label{eqn: nonlinear cartpole system}
	\ddot{x} &= \frac{u + ml\omega^{2} \sin(\theta)}{m_{t}} \\
	&- \frac{ml (g \sin(\theta) - \cos(\theta)) (\frac{u + ml\omega^{2} \sin(\theta)}{m_{t}})}{l(\frac{4}{3} - m \frac{\cos^{2}(\theta)}{m_{t}})} \frac{\cos(\theta)}{m_{t}} \\
	\ddot{\theta} &= \frac{g \sin(\theta) - \cos(\theta) (\frac{u + ml\omega^{2} \sin(\theta)}{m_{t}})}{l(\frac{4}{3} - m \frac{\cos^{2}(\theta)}{m_{t}})}
	\frac{\cos(\theta)}{m_{t}}
\end{split}
\end{align}
where $g = 9.8$ is the gravitational constant, the pole mass is $m = 0.1$, half the pole's length is $l = 0.5$, and $m_{t} = 1.1$ is the total mass.
The dynamics are then discretized in time with sampling time $T_{s} = 0.015$ s.
The control input, $u \in \lbrace -10, 10 \rbrace$, which affects the lateral position of the cart, is chosen by the neural network controller \cite{manzanas2019arch}. 
This means the controller is less ``smooth'' than the neural network controller for the linearized cart-pole system, because when the pendulum is near vertical, the controller rapidly switches between a high positive and negative control input value. Thus, the velocity components of the system state, $\dot{x}$ and $\dot{\theta}$, will not stabilize to zero, as shown in Figure \ref{fig: cartpole nonlinear trajectory}. 

As before, we added a Gaussian disturbance $\boldsymbol{w}_{k} \sim \mathcal{N}(0, \Sigma)$ with $\Sigma = 0.01 I$ to the dynamical equations and represent the system as a Markov control process. The benchmark is defined such that the pole angle $\theta$ will remain within $[-\pi/6, \pi/6]$, while the other variables are unconstrained. 
As such, we define the target set $\mathcal{T}$ such that $\theta \in [-0.05, 0.05]$ as with the linearized cart-pole system, but define the safe set $\mathcal{K}$ to be the entire state space, meaning there is no unsafe region. 
\begin{equation}
    \mathcal{T} = \lbrace x \in \mathbb{R}^{4} \mid \vert x_{3} \vert \leq 0.05 \rbrace
\end{equation}
This allows us to analyze the behavior of the controller to reach a pre-specified objective without enforcing constraints on the system before the terminal time. 
We then simulated $10$ trajectories from initial conditions sampled uniformly from the ranges specified above, and collected a sample $\mathcal{S}$ of $M = 10{,}000$ observations. 
Then, we computed the safety probabilities using Lemma \ref{lem: terminal-hitting time problem} with $N = 600$. The results are shown in Figure \ref{fig: cartpole nonlinear}(a). 
We then computed the finite sample bounds using Theorem \ref{thm: main result} for $\delta = 0.1$, and plotted the lower bound on the safety probabilities in Figure \ref{fig: cartpole nonlinear}(b).  

As expected, we see that the safety probabilities show that the nonlinear pendulum system is able to stabilize a pendulum starting within a small range of $\theta$. Interestingly, it is revealed in Figure \ref{fig: cartpole nonlinear}(a) that the controller is not completely symmetric, meaning it has a higher chance of stabilizing the pendulum for positive values of $\theta$ than negative values of $\theta$.
Fig. \ref{fig: cartpole nonlinear}(b) shows the probabilistic lower bound on safety probabilities computed using Theorem \ref{thm: main result} with $\delta = 0.1$.

\section{Conclusion}
\label{section: conclusion}

We provided state- and input-based finite sample bounds for the stochastic reachability probability constructed via conditional distribution embeddings.  
Our approach is based on an application of statistical learning theory, that relates the observed data to the quality of the approximation of the stochastic reachability probability at a given state and input.  
This approach enables rigorous bounds on model-free stochastic reachability.  
We validated our approach on a nonlinear dynamical system with a neural net controller, and numerically characterized our approach on the stochastic double integrator.  



\bibliographystyle{plain}
\bibliography{bibliography}

\end{document}